\documentclass[12pt]{amsart}
\usepackage{amsmath,amssymb,amstext,amsthm,amscd}

\theoremstyle{plain}
\newtheorem{thm}{Theorem}
\newtheorem{pro}[thm]{Proposition}
\newtheorem{cor}{Corollary}
\newtheorem*{thm2}{Theorem}
\newtheorem*{pro2}{Proposition}
\theoremstyle{definition}
\newtheorem*{ack}{Acknowledgements}

\def\pf{{\em Proof.}\ \,}
\def\su{\subseteq}

\def\fun{\rightarrow}
\def\lfun{{\longrightarrow}}

\def\Th{\Theta}

\def\bC{{\mathbb{C}}}
\def\bZ{{\mathbb{Z}}}

\def\cL{{\mathcal L}}
\def\cO{{\mathcal O}}
\def\SU{{\mathcal S}{\mathcal U}}
\def\cU{{\mathcal U}}

\def\Gr{{\textrm{Gr}}}
\def\rk{{\textrm{rk}}}

\begin{document}

\title{A lower bound for the dimension of the base locus of the generalized
theta divisor}

%    Information for first author
\author{Daniele Arcara}
\address{Department of Mathematics, University of Utah,
155 S. 1400 E., Room 233, Salt Lake City, UT 84112-0090, USA}
\email{arcara@math.utah.edu}

\begin{abstract}
We produce a lower bound for the dimension of the base locus of the
generalized theta divisor $\Theta_r$ on the moduli space $\SU_C(r)$ of
semistable vector bundles of rank $r$ and trivial determinant on a smooth
curve $C$ of genus $g\geq2$.
\end{abstract}

\maketitle

\section{Introduction}

Let $C$ be a smooth irreducible complex projective curve of genus $g\geq2$.
Let $\cU_C(r)$ be the moduli space of (S-equivalence classes of) semi-stable
vector bundles of rank $r$ and degree $0$, and let $\SU_C(r)$ be the moduli
space of (S-equivalence classes of) semi-stable vector bundles of rank $r$ and
trivial determinant.

The Picard group of $\SU_C(r)$ is generated by an ample line bundle (see
[DreNar89]), that we shall denote by $\cL_r$.
A divisor $\Th_r$ on $\SU_C(r)$ such that $\cL_r=\cO_{\SU_C(r)}(\Th_r)$ is
called a generalized theta divisor.
We are interested in the base locus of its linear system.

A vector bundle $E\in\SU_C(r)$ is in the base locus of the generalized theta
divisor if and only if $H^0(E\otimes L)\neq0$ for every line bundle $L$ on $C$
of degree $g-1$ (see [Bea88] for $r=2$ and [BeaNarRam89] in general).
Raynaud studied bundles with a similar property in [Ray82], and Beauville
summarizes his results as follows in [Bea95].

\begin{thm2}[Raynaud]
$(a)$ For $r=2$, the linear system $|\Th_2|$ has no base points.

$(b)$ For $r=3$, $|\Th_3|$ has no base points if $g=2$, or if $g\geq3$ and $C$
is generic.

$(c)$ Let $n$ be an integer $\ge2$ dividing $g$.
For $r=n^g$, the linear system $|\Th_r|$ has base points.
\end{thm2}

For part $(c)$, Raynaud actually constructs, for every $n\ge2$ and $g$, a
vector bundles of rank $n^g$ and slope $g/n$ without the property that he
calls $(*)$.
The bundles without the property $(*)$ can be easily used to produce vector
bundles in the base locus of the generalized theta divisor if their slope is
integral, hence part $(c)$ of the theorem above.
Note also that Popa generalized Raynaud's construction in [Pop02].

Using the dual $E_L$ of the kernel of the evaluation map $e_L$ for a line
bundle $L$ generated by its global sections, and its exterior powers, Popa
[Pop99] and Schneider [Sch03] proved the existence of other vector bundles in
the base locus of the generalized theta divisor.

In particular, Schneider defines a condition $(R)$, which implies Raynaud's
condition $(*)$, as follows:
A vector bundles $E$ has the property $(R)$ if, for every $n\in\bZ$ and any
generic line bundle $L$ of degree $n$, $H^0(E\otimes L)=0$ or
$H^1(E\otimes L)=0$.
He then proves the following proposition.

\begin{pro2}[Schneider]
Let $C$ be a smooth complex projective curve of genus $g\ge2$.
If $L$ is a line bundle of degree greater than or equal to $2g+2$, then
$\Lambda^pE_L$ does not verify $(R)$ for every
$p\in\{2,\dots,\emph{rk}(E_L)-2\}$.
\end{pro2}

Under the assumption of the proposition, $E_L$ is stable (see [EinLaz92]), and
therefore $\Lambda^pE_L$ is semi-stable.
Whenever the slope of $\Lambda^pE_L$ is integral, this easily produces
examples of vector bundles in the base locus of the generalized theta divisor.

As our first result, we prove that every vector bundles without the property
$(R)$ ``produces'' a vector bundle in the base locus of the generalized theta
divisor, hence making it possible to use all of the bundles studied by
Raynaud, Popa, and Schneider, even the ones with non-integral slope.

\begin{thm}\label{slope}
If $E$ is a semi-stable vector bundle of rank $r$ on $C$ which does not
satisfy the property $(R)$, then the base locus of $|\Th_r|$ is non-empty.
\end{thm}

As a corollary, using Raynaud's and Schneider's results, we obtain the
following corollary.

\begin{cor}\label{nonempty}
The base locus of $|\Th_r|$ is non-empty for
$ r = 2^g $ and $ r = (g+1)(g+2)/2. $
\end{cor}

As Popa points out in [Pop99], this implies that the base locus is also
non-empty for any bigger rank (if $E$ is in the base locus of $|\Th_r|$, just
take $E\oplus\cO_C^{\oplus n}$).
If we let $r_0$ be the lowest rank such that the base locus of $|\Th_r|$ is
non-empty, the corollary above can be restated as
$$ r_0 \leq \min \left\{ 2^g, \frac{(g+1)(g+2)}{2} \right\}. $$

We produce the following lower bound for the dimension of the base locus of
$|\Th_r|$ on a curve of genus $g$:

\begin{thm}\label{lowerbound}
Let $C$ be a smooth complex projective curve of genus $g$.
Then the dimension of the base locus of $|\Th_r|$ is at least
$(r-r_0)^2(g-1)+1$, where $r_0$ is the minimum rank for which the base
locus of the generalized theta divisor is non-empty.
\end{thm}

\begin{ack}
I would like to thank E.\ Izadi for suggesting the problem, and V.\ Alexeev
for pointing out a way to simplify the proof of theorem \ref{lowerbound}.
\end{ack}

\section{Proof of theorem \ref{slope}}

Schneider proves in [Sch03] that a vector bundle $E$ satisfies the condition
$(R)$ if and only if it satisfies the two following conditions:

(1) $H^1(E\otimes L)=0$ for a generic line bundle $L$ of degree
$g-1-\lfloor\mu(E)\rfloor$;

(2) $H^0(E\otimes L)=0$ for a generic line bundle $L$ of degree
$g-1-\lceil\mu(E)\rceil$.

Let $E$ be a semi-stable vector bundle of rank $r$ which does not satisfy the
property $(R)$.
If $\mu(E)$ is an integer, then $E\otimes L$ is in the base locus of
$|\Th_r|$, where $L$ is a line bundle of degree $g-1-\mu(E)$ such that
$L^{-r}\simeq\det E$.
If $\mu(E)$ is not an integer, there are two cases.

Case I:
For every line bundle $L$ of degree $g-1-\lfloor\mu(E)\rfloor$,
$H^1(E\otimes L)\neq0$.

Let $\bC_p$ be a skyscreaper sheaf of degree $1$ supported at a point $p$ of
$C$, let $E\fun\bC_p$ be a non-zero map, and let $E'$ be the kernel:
$ 0 \fun E' \fun E \fun \bC_p \fun 0. $
Since $H^1(\bC_p)=0$, $E'$ also satisfies the condition that
$H^1(E'\otimes L)\neq0$ for every line bundle $L$ of degree
$g-1-\lfloor\mu(E)\rfloor$.
Moreover, $\lfloor\mu(E')\rfloor=\lfloor\mu(E)\rfloor$.
There are now two subcases.

Subcase I.1:
$E'$ is semi-stable.
Then $E'$ is a semi-stable vector bundle of slope $\mu(E')<\mu(E)$ with
$\lfloor\mu(E')\rfloor=\lfloor\mu(E)\rfloor$ such that
$H^1(E'\otimes L)\neq0$ for every line bundle $L$ of degree
$g-1-\lfloor\mu(E')\rfloor$.

Subcase I.2:
$E'$ is not semi-stable.
Let $\mu$ be the maximum slope of a vector sub-bundle of $E'$, and let $F'$ be
a sub-bundle of maximal rank among all of the sub-bundles of slope $\mu$.
Then there exists a short exact sequence $0\fun F'\fun E'\fun G'\fun0$ with
$F'$ and $G'$ stable vector bundles.
By semi-continuity, we obtain that either $H^1(F'\otimes L)\neq0$ for every
line bundle $L$ of degree $g-1-\lfloor\mu(E)\rfloor$ or
$H^1(G'\otimes L)\neq0$ for every line bundle $L$ of degree
$g-1-\lfloor\mu(E)\rfloor$.
Let us show that $\mu(F')$ and $\mu(G')$ are both $\geq\lfloor\mu(E)\rfloor$.
Clearly, $\mu(F')>\mu(E')\geq\lfloor\mu(E)\rfloor$.
For $G'$, note that it is contained in $G=E/F'$, and
$\mu(G)\geq\mu(E)>\lfloor\mu(E)\rfloor$.
Therefore, $\mu(G')=\mu(G)-1/\rk(G)\geq\lfloor\mu(E)\rfloor$.
Since $\mu(F')$ and $\mu(G')$ are clearly $\leq\mu(E)$, this shows that
$\lfloor\mu(F')\rfloor=\lfloor\mu(G')\rfloor=\lfloor\mu(E)\rfloor$.
Therefore, either $F'$ or $G'$ is a semi-stable vector bundle $E''$ of slope
$\mu(E'')\geq\lfloor\mu(E)\rfloor$ and rank $\rk(E'')<r$ with
$\lfloor\mu(E'')\rfloor=\lfloor\mu(E)\rfloor$ such that
$H^1(E''\otimes L)\neq0$ for every line bundle $L$ of degree
$g-1-\lfloor\mu(E'')\rfloor$.

We can now continue our process by replacing $E$ with either $E'$, if it is
semi-stable, or with the $E''$ constructed in the case when $E'$ is not
semi-stable.
The process eventually ends when the slope becomes integral.
This happens because the slopes of the vector bundles constructed at each step
is bounded below by $\lfloor\mu(E)\rfloor$ and at each step either the slope
or the rank is decreasing.

Case II:
For every line bundle $L$ of degree $g-1-\lceil\mu(E)\rceil$,
$H^0(E\otimes L)\neq0$.

Let $\bC_p$ be a skyscreaper sheaf of degree $1$ supported at a point $p$ of
$C$, let $E'$ be a non-trivial extension of $\bC_p$ by $E$:
$ 0 \fun E \fun E' \fun \bC_p \fun 0. $
It is clear that $H^0(E'\otimes L)\neq0$ for every line bundle $L$ of degree
$g-1-\lceil\mu(E)\rceil$.
Moreover, $\lceil\mu(E')\rceil=\lceil\mu(E)\rceil$.
As before, there are two subcases.

Subcase II.1:
$E'$ is semi-stable.
Then $E'$ is a semi-stable vector bundle of slope $\mu(E')>\mu(E)$ with
$\lceil\mu(E')\rceil=\lceil\mu(E)\rceil$ such that
$H^0(E'\otimes L)\neq0$ for every line bundle $L$ of degree
$g-1-\lceil\mu(E')\rceil$.

Subcase II.2:
$E'$ is not semi-stable.
Construct $F'$ and $G'$ as in subcase I.2.
By semi-continuity, we obtain that either $H^0(F'\otimes L)\neq0$ for every
line bundle $L$ of degree $g-1-\lceil\mu(E)\rceil$ or
$H^0(G'\otimes L)\neq0$ for every line bundle $L$ of degree
$g-1-\lceil\mu(E)\rceil$.
Let us show that $\mu(F')$ and $\mu(G')$ are both $\leq\lceil\mu(E)\rceil$.
Clearly, $\mu(G')<\mu(E')\leq\lceil\mu(E)\rceil$.
Let $F=E\cap F'$.
Since $\mu(F')>\mu(E')>\mu(E)$, and $E$ is semi-stable, $F\neq F'$.
Since $E$ is semi-stable, $\mu(F)\leq\mu(E)<\lceil\mu(E)\rceil$, and
therefore $\mu(F')=\mu(F)+1/\rk(F)\leq\lceil\mu(E)\rceil$.
Therefore, either $F'$ or $G'$ is a semi-stable vector bundle $E''$ of slope
$\mu(E'')\leq\lceil\mu(E)\rceil$ and rank $\rk(E'')<r$ with
$\lceil\mu(E'')\rceil=\lceil\mu(E)\rceil$ such that
$H^1(E''\otimes L)\neq0$ for every line bundle $L$ of degree
$g-1-\lceil\mu(E'')\rceil$.

We can now continue our process by replacing $E$ with either $E'$, if it is
semi-stable, or with the $E''$ constructed in the case when $E'$ is not
semi-stable.
The process eventually ends when the slope becomes integral.
This happens because the slopes of the vector bundles constructed at each step
is bounded above by $\lceil\mu(E)\rceil$ and at each step either the slope is
increasing or the rank is decreasing.

To conclude the proof of the theorem, note that, if the vector bundle in the
base locus of the generalized theta divisor constructed has rank $r'<r$, we
can always produce one in the base locus of the generalized theta divisor by
doing a direct sum with $r-r'$ copies of $\cO_C$.

\section{Proof of theorem \ref{lowerbound}}

To simplify the proof of theorem \ref{lowerbound}, let us point out the
following result, which is probably well-known.
We prove it here because we could not find a proof in the literature.

\begin{pro}\label{prosta}
Every vector bundle in the base locus of $|\Th_{r_0}|$ is stable.
\end{pro}

\pf
Let $E$ be a vector bundle in the base locus of $|\Theta_{r_0}|$.
Then $E$ is in the same equivalence class as its associated grading
$\oplus_{i=1}^k\Gr_i$ from its Jordan-H\"older filtration.
Since $H^0(E\otimes L)\ne0$ if and only if there exists an $i$ such that
$H^0(\Gr_i\otimes L)\ne0$, by semicontinuity there exists an $i$ such that
$\Gr_i$ is in the base locus of $|\Theta_{\rk\Gr_i}|$.
By the minimality of $r_0$, $\rk\Gr_i=r_0$, $k=1$, and $E=\Gr_1$ is
stable.
\qed

We are now ready to prove the theorem.
Let $E$ be a vector bundle in the base locus of $|\Theta_{r_0}|$.
Then $E$ is stable by proposition \ref{prosta}.
Let $n=r-r_0$, and consider the morphism
$$ \varphi \colon JC \times \cU_C(n) \lfun \cU_C(r) \quad
\varphi(L, F) = (E \otimes L) \oplus F, $$
where $JC$ is the Jacobian of $C$.
There is a natural embedding $\SU_C(r)\su\cU_C(r)$, and we claim that, if we
let
$ B = \varphi(JC \times \cU_C(n)) \cap \SU_C(r), $
then $\dim B=n^2(g-1)+1$, and $B$ is contained in the base locus of $|\Th_r|$.

Let $E'\in B$.
Then, by definition of $B$, $E'$ contains $E\otimes L$ for some $L\in JC$.
Then, for every line bundle $L'$ of degree $g-1$, $H^0(E'\otimes L')$ contains
$H^0(E\otimes(L\otimes L'))$, which is non-zero since $E$ is in the base locus
of the generalized theta divisor.
Therefore, $E'$ is in the base locus of $|\Th_r|$.

To prove the claim about the dimension of $B$, let $A$ be the preimage of
$\SU_C(r)$ under $\varphi$, i.e.,
$A = \{ (L, F) \in JC \times \cU_C(n) \mid L^r \otimes \det F \simeq \cO_C\}.$
Then $B=\varphi(A)$.

To finish the proof, it suffices to prove that
$\dim A=\dim\cU_C(n)=n^2(g-1)+1$ and that the general fiber of $\varphi|_A$ is
finite.

Consider the second projection $A\fun\cU_C(n)$: it is surjective because for
every $F\in\cU_C(n)$ there exists $L\in JC$ such that
$L^r\otimes\det F\simeq\cO_C$.
Moreover, there are only a finite number of such $L$'s, and this proves the
claim about the dimension of $A$.

Let $E'$ be a generic element of $B$.
It is of the form $\varphi(L,F)=(E\otimes L)\oplus F$ for some $(L,F)\in A$,
and we can assume in what follows that $F$ is stable, since a generic element
of $\cU_C(n)$ is stable.
Then the associated grading for $E'$ is $(E\otimes L)\oplus F$ itself, being
both $E\otimes L$ and $F$ stable.

Let us show that there are only a finite number of $(L',F')\in A$ such that
$E'$ is S-equivalent to $(E\otimes L')\oplus F'$.
Indeed, if this is the case, then $F'$ must be stable, and there are only two
possibilities:
(I) $F\simeq F'$ and $E\otimes L\simeq E\otimes L'$.
Then $F'$ is unique in $\cU_C(n)$, and $L^r\simeq(L')^r$ implies that there
are only finitely many such $L'$'s.
(II) $E\otimes L\simeq F'$ and $F\simeq E\otimes L'$.
Then again $F'$ is uniquely determined, and $(L')^r\simeq\det F$ implies that
there are only finitely many such $L'$'s.

\end{document}